\documentclass[12pt]{amsart}
\usepackage{amscd,amsmath,amssymb,amsfonts}
\usepackage[cmtip, all]{xy}
\begin{document}
\newtheorem{thm}[subsection]{Theorem}
\newtheorem{prop}[subsection]{Proposition}
\newtheorem{lem}[subsection]{Lemma}
\newtheorem{cor}[subsection]{Corollary}
\theoremstyle{definition}
\newtheorem{Def}[subsection]{Definition}
\theoremstyle{remark}
\newtheorem{rem}[subsection]{Remark}
\newtheorem{rems}[subsection]{Remarks}
\newtheorem{ex}[subsection]{Example}
\newtheorem{exs}[subsection]{Examples}
\numberwithin{equation}{section}
\newcommand{\an}{{\rm an}}
\newcommand{\alg}{{\rm alg}}
\newcommand{\cl}{{\rm cl}}
\newcommand{\Alb}{{\rm Alb}}
\newcommand{\CH}{{\rm CH}}
\newcommand{\mc}{\mathcal}
\newcommand{\mb}{\mathbb}
\newcommand{\surj}{\twoheadrightarrow}
\newcommand{\inj}{\hookrightarrow}
\newcommand{\red}{{\rm red}}
\newcommand{\codim}{{\rm codim}}
\newcommand{\rank}{{\rm rank}}
\newcommand{\Pic}{{\rm Pic}}
\newcommand{\Div}{{\rm Div}}
\newcommand{\Hom}{{\rm Hom}}
\newcommand{\im}{{\rm im}}
\newcommand{\Spec}{{\rm Spec \,}}
\newcommand{\sing}{{\rm sing}}
\newcommand{\Char}{{\rm char}}
\newcommand{\Tr}{{\rm Tr}}
\newcommand{\Gal}{{\rm Gal}}
\newcommand{\Min}{{\rm Min \ }}
\newcommand{\Max}{{\rm Max \ }}
\newcommand{\supp}{{\rm supp}\,}
\newcommand{\0}{\emptyset}
\newcommand{\sA}{{\mathcal A}}
\newcommand{\sB}{{\mathcal B}}
\newcommand{\sC}{{\mathcal C}}
\newcommand{\sD}{{\mathcal D}}
\newcommand{\sE}{{\mathcal E}}
\newcommand{\sF}{{\mathcal F}}
\newcommand{\sG}{{\mathcal G}}
\newcommand{\sH}{{\mathcal H}}
\newcommand{\sI}{{\mathcal I}}
\newcommand{\sJ}{{\mathcal J}}
\newcommand{\sK}{{\mathcal K}}
\newcommand{\sL}{{\mathcal L}}
\newcommand{\sM}{{\mathcal M}}
\newcommand{\sN}{{\mathcal N}}
\newcommand{\sO}{{\mathcal O}}
\newcommand{\sP}{{\mathcal P}}
\newcommand{\sQ}{{\mathcal Q}}
\newcommand{\sR}{{\mathcal R}}
\newcommand{\sS}{{\mathcal S}}
\newcommand{\sT}{{\mathcal T}}
\newcommand{\sU}{{\mathcal U}}
\newcommand{\sV}{{\mathcal V}}
\newcommand{\sW}{{\mathcal W}}
\newcommand{\sX}{{\mathcal X}}
\newcommand{\sY}{{\mathcal Y}}
\newcommand{\sZ}{{\mathcal Z}}
\newcommand{\A}{{\mathbb A}}
\newcommand{\B}{{\mathbb B}}
\newcommand{\C}{{\mathbb C}}
\newcommand{\D}{{\mathbb D}}
\newcommand{\E}{{\mathbb E}}
\newcommand{\F}{{\mathbb F}}
\newcommand{\G}{{\mathbb G}}
\renewcommand{\H}{{\mathbb H}}
\newcommand{\I}{{\mathbb I}}
\newcommand{\J}{{\mathbb J}}
\newcommand{\M}{{\mathbb M}}
\newcommand{\N}{{\mathbb N}}
\renewcommand{\P}{{\mathbb P}}
\newcommand{\Q}{{\mathbb Q}}
\newcommand{\R}{{\mathbb R}}
\newcommand{\T}{{\mathbb T}}
\newcommand{\U}{{\mathbb U}}
\newcommand{\V}{{\mathbb V}}
\newcommand{\W}{{\mathbb W}}
\newcommand{\X}{{\mathbb X}}
\newcommand{\Y}{{\mathbb Y}}
\newcommand{\Z}{{\mathbb Z}}

\newcommand{\Nm}{{\operatorname{Nm}}}
\newcommand{\NS}{{\operatorname{NS}}}
\newcommand{\id}{{\operatorname{id}}}
\newcommand{\Zar}{{\text{\rm Zar}}} 

\newcommand{\zdz}{z\nabla_{\partial_z}}
\newcommand{\res}{{\rm res}}
\newcommand{\Rea}{{\rm Re\,}}

\addtocounter{section}{-1}
\title[Semistable bundles and reducible representations]
{ Semistable bundles on curves and reducible
representations of the fundamental group}

\author{H\'el\`ene Esnault
\and Claus Hertling}

\address{
Universit\"{a}t GH Essen\\
FB6 Mathematik und Informatik\\
45117 Essen\\
Germany\\
\and
Mathematisches Institut der Universit\"{a}t Bonn\\
Beringstra{\ss}e 1\\
53115 Bonn \\
Germany}

\thanks{Partially supported by the
DFG-Forschergruppe ``Arithmetik und Geometrie''}

\email{ esnault@uni-essen.de, hertling@math.uni-bonn.de}




\maketitle

\section{Introduction}

A.A. Bolibruch \cite{Bo} and V.P. Kostov \cite{K} showed
independently that if $\rho:
\pi_1(\P^1 \setminus \Sigma) \to GL(n,\C)$ is an irreducible
representation of the fundamental group, then there is an
algebraic bundle $E$ together with an algebraic connection 
$\nabla: E \to \Omega^1(\log \Sigma)\otimes E$ with underlying
local system $\rho$, with the property that $E \cong \oplus_1^n
\sO$ is algebraically trivial.
 Equivalently, $E$ can be taken to be the twist $L\otimes
(\oplus_1^n \sO)$ of a line bundle $L$
by an algebraically trivial bundle. 
Those twists are the unique semistable bundles on $\P^1$. In 
\cite{EV}, it is indeed proven that if $\P^1$ is replaced by a
smooth projective complex curve of higher genus, the theorem remains
true in this form: there is a $(E,\nabla)$ as above with $E$
semistable of degree 0 (and also with $E$ semistable of any
degree, even if it is not emphasized in the article).
Let us call here for short such an $(E,\nabla)$ a {\it realization}
of $\rho$. Note, it is crucial to require $\nabla$ to have poles only
along $\Sigma$. If one allows one more pole, then one can for
example on $\P^1$ trivialize $E$ even with parameters (see
\cite{BE}, section 4).

On the other hand, on $\P^1_\C$,
A. Bolibruch \cite{AB}  constructed representations which cannot
be realized on the trivial bundle. The purpose of this note is
to show that on a higher genus Riemann surface, there are
representations $\rho$ which cannot be realized on a semistable
bundle. We show:

\begin{thm} \label{mainthm} Let $X$ be a smooth projective complex curve
 of genus $g$ and let 
$\Sigma \subset X$ be a finite nonempty set.

If $(g=0,|\Sigma|\geq 3,n\geq 4)$ or $(g\geq 1,|\Sigma|\geq 1,n\geq 5)$
there exists a representation
$$
\rho :\pi_1 (X\setminus \Sigma) \to GL(n,\C)
$$
which cannot be realized by an algebraic connection
$$
\nabla : E \to \Omega^1_{X} (\log \Sigma) 
                            \otimes E
$$
with logarithmic poles along $\Sigma$ and with $E$ semistable.
\end{thm}

The proof is an adaptation of Bolibruch's ideas to the higher
genus case, together with the use of Gabber's algebraic 
view (\cite{EV}, section 1) on the Bolibruch-Kostov theorem.

Finally, let us remark that Deligne  extensions $(E,\nabla)$
(\cite{Del}) are very natural in geometry. They are not
compatible with pull-backs, but appear as direct
images of connections, for example, Gau{\ss}-Manin connections of
semistable families are Deligne extensions with nilpotent
residues. On the other hand, the Gau{\ss}-Manin bundles tend to
be highly instable, as they contain a positive Hodge subbundle.
Thus it is not clear what is the r\^{o}le of semistability for
realization of monodromy (see one computation in section 5). 

\section{Bolibruch's construction}

Throughout the note we use the following notations.
\begin{list}{}{}
\item[(1.A)] $X$ is a smooth projective complex
curve, $\Sigma =\{p_1,...,p_m\}\subset X$
is finite and nonempty.
$$
\rho :\pi_1 (X\setminus \Sigma) \to GL(n,\C)
$$
is a representation. $E$ is a vector bundle on $X$
of rank $n$,
$$
\nabla : E \to \Omega^1_{X} (\log \Sigma) 
                            \otimes E
$$
is a logarithmic connection on $E$ with underlying local system
$\rho$. We call $(E,\nabla)$ a {\it realization} of $\rho$.
 
The eigenvalues of the residue endomorphism
$$
\res_{p_i}(\nabla): E\otimes \C (p_i) \to E\otimes \C (p_i)
$$
at $p_i$ are called $\beta_{i1},...,\beta_{in}$. They are ordered such that
$$
\Rea \beta_{i1}\leq ... \leq \Rea \beta_{in}\ .
$$
\end{list}

The following theorem is the key to Bolibruch's examples.

\begin{thm} \label{ss}
Let $X, \Sigma, \rho, E, \nabla$, and $\beta_{ij}$ be as in (1.A).
Suppose that $E$ is semistable, that $\rho$ is
reducible, and that for each $i\in \{1,...,m\}$ the monodromy of $\rho$
around $p_i$ has only one eigenvalue $\lambda_i$ and only one Jordan block.

Then $\beta_{i1}=...=\beta_{in}=:\beta_i$ for all $i$ and the slope
$\mu(E)=\frac{\deg (E)}{\rank (E)}$ satisfies
$$
e^{2\pi \sqrt{-1}\mu(E)} = \prod_i \lambda_i\ .
$$
Said differently, $(E,\nabla)$ is the Deligne extension
characterized by the property that
$(\res_{p_i}(\nabla) -\beta_i I)$ nilpotent. 
\end{thm}

We first prove a local statement.
\begin{lem} \label{lem1:2}
Let $j: U=X\setminus \{p\}\inj X$ be the embedding 
of the complement of a point on a smooth analytic contractible curve.
Let $(E,\nabla)$ be a regular
connection on $U$, such that the underlying local
monodromy has only one eigenvalue and one Jordan block.
Thus $E$ has a filtration $E_i\subset
E_{i+1}$ stabilized by $\nabla$.
Let $F\subset j_*E$ be a bundle such that
$\nabla|_{F}$ has logarithmic poles in $\{p\}$, and let us
denote by $\beta_\ell$ its eigenvalues, ordered such that $(\beta_{\ell
+1}-\beta_\ell) \in \N$. Let
$F_i:=j_*E_i \cap F$. Then $F_i\subset F$ is a subbundle,
$\nabla|_{F_i}$ has logarithmic poles
and its residues are precisely $\{\beta_1,\ldots,\beta_i\}$.
\end{lem}

\begin{proof} 
If the rank of $E$ is 1, there is of course nothing to prove. 
Since $X$ is smooth and has dimension 1, $F_i\subset F$ is a
subbundle. As is well known, $\nabla|_F$ stabilizes $F_i$, for it takes
values in $\Omega^1(\log \{p\})\otimes E \cap \Omega^1(\log
\{p\})\otimes j_*E_i$. Furthermore, each $(F_i, E_i)$ satisfies
the same assumptions as $(F,E)$. 
Let us thus first consider $F_2$. If $(\beta_1-\beta_2) \in
\N\setminus\{0\}$, one performs
Gabber's construction ( \cite{EV}, section 1): $F_i$
embedds into $F'_i$, for $i=1,2$, with cokernel $\C(p)$, such
that $\nabla_{F}$ extends as a logarithmic connection, 
the residue of which has the new eigenvalues  $\beta_1-1, \beta_2$.
Furthermore if
the local generators in $\{p\}$ of $F$ are $e_1, e_2$ with
$e_1$ generating $F_1$ and $e_2 \otimes \C(p)$ being an
eigenvector to $\beta_2$, then the new generators are $
(\frac{e_1}{z}, e_2)$. In this basis, one has
 ${\res}'_{p}-\res_{p}= {\rm
diag}(-1,0)$. Repeating the procedure $(\beta_1-\beta_2)$ times,
one reaches a new $F_2$ with residue $={\rm diag} (\beta_2,
\beta_2)$. Now by Deligne \cite{Del} again, this implies that
the local monodromy underlying  $E_2$ is diagonal (actually even
a homothety), a
contradiction. Thus $(\beta_1-\beta_2) \le 0$. Replacing now $E$
by $E/E_1$, one proceeds inductively.
\end{proof}

\bigskip
{\it Proof of theorem \ref{ss}.}
Because $\rho$ is reducible, the local system 
$\ker(\nabla|_{X\setminus \Sigma})$
contains a local subsystem $V$ of some rank $\ell$ with $0<\ell<n$.
Let $(\sV, \nabla|_{\sV}) \subset 
(E|_{X\setminus \Sigma},\nabla|_{X\setminus \Sigma})$ be the 
induced algebraic regular connection. Let 
$j:X\setminus \Sigma\to X$ be the inclusion and $$
F:= j_{*}(\sV)\cap E \subset E\ .
$$
By lemma \ref{lem1:2}, $F$ is a subbundle of $E$ and $\nabla$ restricts to a 
logarithmic connection on $F$ with residue eigenvalues 
$\beta_{i1},...,\beta_{i\ell}$ at $p_i$.

On the other hand, one has (see \cite{EV1}, appendix B, for example)
$$\sum_{i=1}^m \sum_{j=1}^n \beta_{ij}=-{\rm deg}(E).$$
Thus  semistability of $E$ implies
$$
\frac{1}{k}\sum_{i=1}^m\sum_{j=1}^\ell \beta_{ij} = -\mu(F) \geq -\mu(E) =
\frac{1}{n}\sum_{i=1}^m\sum_{j=1}^n \beta_{ij}\ .
$$
Together with $\beta_{ij}- \beta_{i\ell} \ge 0 \ \text{for} \
j\ge \ell$
this shows
$$
\beta_{i1}=...=\beta_{in}
$$
and
$$
e^{2\pi \sqrt{-1} \mu(E)} = \prod_{i=1}^me^{-2\pi \sqrt{-1}\beta_{i1}}
=\prod_{i=1}^m \lambda_i \ .
$$
\hfill \qed

\section{Examples of reducible representations}

Here we list several representations as in theorem \ref{ss}.

As in 1.A, 
$X$ is a smooth projective complex curve of genus $g$, 
and $\Sigma =\{p_1,...,p_m\}\subset X$ is a finite nonempty set.
One chooses a system of paths $a_1,b_1,...,a_g,b_g$ and $c_1,...,c_m$ with
a common base point such that $\pi_1(X\setminus \Sigma)$ is generated by them with
the single relation
$$
a_1b_1a_1^{-1}b_1^{-1}\cdot ...\cdot a_gb_ga_g^{-1}b_g^{-1}\cdot c_1...c_m=0\ .
$$
Here $c_i$ is a loop around $p_i$.
Then a representation $\rho:\pi_1(X\setminus \Sigma)\to GL(n,\C)$ is given by
$n\times n$-matrices $A_1,B_1,...,A_g,B_g$ and $C_1,...,C_m$
which satisfy the same relation. 

The starting point is to overtake Bolibruch's examples by
setting $A_i=B_i={\rm id}$, ((2.1)--(2.3)) and then to
modify ((2.4)-(2.5)).

The following representations are all reducible, and the local monodromies
$C_i$ around the points $p_i$ have only one eigenvalue $\lambda_i$ 
and only one Jordan block.

\medskip
\noindent
{\bf (2.1)} $\rho^{(1)}$: Choose $\nu_i\in \C \setminus\{0\}$, 
$i=1,...,m$, with 
$\sum_{i=1}^m\nu_i=0$ and a nilpotent $n\times n$-matrix 
$N^{(1)}$ with $\rank\, N^{(1)}=n-1$. Define
$$
C_i^{(1)}:= \exp (\nu_i N^{(1)})\ ,\ \ A_i^{(1)}:=B_i^{(1)}:= \id\ .
$$
Then $\lambda_i=1$.

\medskip
\noindent
{\bf (2.2)} (Bolibruch \cite{AB} Example 5.3.1) 
$\rho^{(2)}$: $m=3,\ n=4, \ A_i^{(2)}:=B_i^{(2)}:=\id$,
\begin{eqnarray*}
C_1^{(2)} := \left( \begin{array}{cccc}
1 & 1 & 0 & 0 \\
0 & 1 & 1 & 0 \\
0 & 0 & 1 & 1 \\
0 & 0 & 0 & 1 
\end{array} \right)\ , \ 
C_2^{(2)} := \left( \begin{array}{cccc}
3 & 1 & 1 & -1 \\
-4 & -1 & 1 & 2 \\
0 & 0 & 3 & 1 \\
0 & 0 & -4 & -1  
\end{array} \right)\ , 
\end{eqnarray*}
\begin{eqnarray*} 
C_3^{(2)} := \left( \begin{array}{cccc}
-1 & 0 & 2 & -1 \\
4 & -1 & 0 & 1 \\
0 & 0 & -1 & 0 \\
0 & 0 & 4 & -1 
\end{array} \right)\ .
\end{eqnarray*}
Then $\lambda_1=\lambda_2=1,\ \lambda_3=-1$.
A semistable bundle $E$ with a logarithmic connection (with poles only in
$\Sigma$) which realizes $\rho^{(2)}$ must have slope
$\mu(E)\equiv \frac{1}{2}\mod \Z$ by theorem \ref{mainthm}.

In the case $g=0$ this is impossible as any semistable bundle
has slope in $\Z$. In that case $\rho^{(2)}$ cannot be realized by a 
logarithmic connection on a semistable bundle (with poles only in $\Sigma$)
and in particular not by a Fuchsian differential system.

\medskip
\noindent
{\bf (2.3)} $\rho^{(3)}$: $m\geq 3, \ n=4,\ A_i^{(3)}:=B_i^{(3)}:=\id$.
Define $N^{(3)}:= \log C_1^{(2)}$,
\begin{eqnarray*}
&& C_1^{(3)}:=...:=C_{m-2}^{(3)}:=\exp (2\pi \sqrt{-1}\frac{1}{2m-4})
\exp(\frac{1}{m-2} N^{(3)})\\
&& C_{m-1}^{(3)}:= -C_2^{2}\ ,\ C_m^{(3)}:=C_3^{(2)}\ .
\end{eqnarray*}
Then $\lambda_1=...=\lambda_{m-2}=\exp (2\pi \sqrt{-1}\frac{1}{2m-4}),\ 
\lambda_{m-1}=\lambda_m=-1, \ \prod_i\lambda_i=-1$.

\medskip
\noindent
{\bf (2.4)} $\rho^{(4)}$: $g\geq 1,\ n=2$. 
$A_i^{(4)}:=B_i^{(4)}:=\id$ for $i\geq 2$. Define
\begin{eqnarray*}
A_1^{(4)}:= \left(\begin{array}{cc} 1&2\\0&1\end{array}\right)\ ,\ 
B_1^{(4)}:= \left(\begin{array}{cc} 1&0\\0&2\end{array}\right)\ ,\ 
C_i^{(4)}:= \left(\begin{array}{cc} 1&\frac{-1}{m}\\0&1\end{array}\right)\ . 
\end{eqnarray*}
Then $\lambda_i=1$.

\medskip
\noindent
{\bf (2.5)} $\rho^{(5)}$: $g\geq 1,\ n$ even, $n\geq 4$. 
$A_i^{(5)}:=B_i^{(5)}:=\id$ for $i\geq 2$. Define
\begin{eqnarray*}
&&\alpha_1:= \left(\begin{array}{cc} 1&2\\0&1 \end{array}\right)\ ,\ 
  \alpha_2:= \left(\begin{array}{cc} 0&1\\0&0 \end{array}\right)\ ,\ 
  \beta   := \left(\begin{array}{cc} 0&1\\1&0 \end{array}\right)\ ,\\
&&\delta_1:= \left(\begin{array}{cc} -3&2\\-2&1 \end{array}\right)\ ,\ 
  \delta_2:= \left(\begin{array}{cc} -4&1\\-1&0 \end{array}\right)\ ,
\end{eqnarray*}
\begin{eqnarray*}
&& A_1^{(5)} := \left( \begin{array}{cccc}
\alpha_1 & \alpha_2 &        & 0       \\
         & \alpha_1 & \ddots &         \\
         &          & \ddots & \alpha_2\\
0        &          &        & \alpha_1
\end{array}\right)\ , \ 
   B_1^{(5)} := \left( \begin{array}{cccc}
\beta    &          &        & 0       \\
         & \beta    &        &         \\
         &          & \ddots &         \\
0        &          &        & \beta 
\end{array}\right)\ . 
\end{eqnarray*}
The matrix
\begin{eqnarray*}
&& A_1^{(5)}B_1^{(5)}{A_1^{(5)}}^{-1}{B_1^{(5)}}^{-1}= 
\left( \begin{array}{cccc}
\delta_1 & \delta_2 &        & {*}     \\
         & \delta_1 & \ddots &         \\
         &          & \ddots & \delta_2\\
0        &          &        & \delta_1
\end{array}\right)
\end{eqnarray*}
has one $n\times n$ Jordan block with eigenvalue $-1$. Define
$$
C_i:=\exp(2\pi \sqrt{-1}\frac{1}{2m})
\exp\left( \frac{-1}{m} 
\log (-A_1^{(5)}B_1^{(5)}{A_1^{(5)}}^{-1}{B_1^{(5)}}^{-1})\right)\ .
$$
Then $\lambda_i=\exp(2\pi \sqrt{-1}\frac{1}{2m})\ , \ \prod_i \lambda_i=-1$.

\section{The proof of theorem \ref{mainthm}}

\begin{thm} \label{sums}
Let $X$ be a smooth projective complex curve, 
$$\Sigma  =   \{p_1,...,p_m\} \subset  X$$ be a 
finite nonempty set, and
$$
\rho_\ell:\pi_1(X\setminus \Sigma)\to GL(n_\ell,\C)\ , \ \ \ell=1,2,
$$
be two representations with the following properties. 
Both representations are reducible.
Each local monodromy has a single eigenvalue and a single Jordan block.
Let $\lambda_i^\ell$ be those eigenvalues  in $p_i$.
Then $\lambda_i^1\neq \lambda_i^2$ and 
$\prod_i \lambda_i^1 \neq \prod_i \lambda_i^2$.

Then $\rho_1\oplus \rho_2$ cannot be realized by a semistable bundle on $X$
with a logarithmic connection with poles only in $\Sigma$.
\end{thm}

\begin{proof}
Let $(E,\nabla)$ be the algebraic regular connection 
on $X\setminus \Sigma$ 
with underlying  $\rho$. Then $(E=E_1 \oplus E_2, \nabla=\nabla_1\oplus
\nabla_2)$, where $\rho_i$ underlies $(E_i,\nabla_i)$. 
Let $F\subset j_*E$ be a bundle such that
$\nabla|_{F}$ has logarithmic poles in $\Sigma$. Then
$F_\ell=j_*E_\ell\cap F\subset F$ is a subbundle, stabilized
by $\nabla|_{F}$.  Let us denote by $\nabla_{F_\ell}$ the
induced connection. Then its residue is the restriction of the
residue of $\nabla_{F}$ to $F_\ell$. Since $\lambda_i^1\neq
\lambda_i^2$ in all points $p_i$, a fortiori none of the
eigenvalues of ${\rm res}_{p_i}(\nabla_{F_1})$ can be an
eigenvalue of ${\rm res}_{p_i}(\nabla_{F_2})$. Consequently, one has
\begin{gather}\label{res}
 {\rm res}_{p_i}(\nabla_{F})= {\rm res}_{p_i}(\nabla_{F_1})
\oplus {\rm res}_{p_i}(\nabla_{F_2}).\end{gather}
On the other hand, 
since $F_1\cap
F_2\subset F$ is torsion free and supported in $\Sigma$, one has
$F_1\cap F_2=0$, thus $F_1\oplus F_2\subset F$ is a locally free
subsheaf, isomorphic to $F$ away of $\Sigma$, and thus isomorphic
to $F$ by the condition (\ref{res}).

If now moreover $F$ is semistable, then $F_\ell$ is semistable as
well, and one has
$\mu(F_1)=\mu(F_2)$. This contradicts theorem \ref{ss}.
\end{proof}

For the proof of theorem \ref{mainthm} one applies theorems
\ref{ss}  and \ref{sums} to 
several combinations of the representations in section 2.\\
$g=0, |\Sigma|\geq 3, n\geq 4$: $\rho^{(1)}$ for $n^{(1)}=n-4$ and 
$\rho^{(3)}$. \\
$g\geq 1, |\Sigma|\geq 1, n $ odd, $n\geq 5$: $\rho^{(1)}$ for $n^{(1)}=1$ and
$\rho^{(5)}$ for $n^{(5)}=n-1$. \\
$g\geq 1, |\Sigma|\geq 1, n $ even, $n\geq 6$: $\rho^{(4)}$ and
$\rho^{(5)}$ for $n^{(5)}=n-2$.

\section{Two-dimensional representations}

W. Dekkers \cite{Dek} showed that, for $X=\P^1_\C$ and a finite nonempty
subset $\Sigma\subset X$, any two-dimensional representation
$\rho:\pi_1(X\setminus \Sigma)\to GL(2,\C)$ can be realized 
on the trivial bundle
with poles only in $\Sigma$.
A.A. Bolibruch gave a simpler proof, using the analogous 
result for  irreducible connections (\cite{Bo}, \cite{K}, \cite{EV}).
We adapt now this to higher genus.

\begin{thm}
Let $X$ be a smooth projective complex curve, $\Sigma\subset X$ a finite 
nonempty set, and 
$$
\rho:\pi_1(X\setminus \Sigma)\to GL(2,\C)
$$
be a two-dimensional representation.

There exists a semistable bundle $E$ of even degree 
with a logarithmic connection
with poles only in $\Sigma$ which realizes $\rho$, but not
necessarily of odd degree.
\end{thm}

\begin{proof}
For $\rho$ irreducible see \cite{EV}. Suppose that $\rho$ is reducible.
Let $(E,\nabla)$ be a vector bundle on $X$ with logarithmic connection
with poles only in $\Sigma$ which realizes $\rho$.

Let $V \subset \ker(\nabla|_{X\setminus \Sigma})$ be a subsystem of rank 1.
We denote by 
$j:X\setminus \Sigma \inj X$ the inclusion and define
$$
F:= j_{*}(V \otimes \sO_{X\setminus \Sigma})\subset E\ .
$$
Then
$$
0\to F \to E \to E/F \to 0
$$
is an exact sequence of bundles, and $F$ and $E/F$ are equipped
with the induced connection. Then $E$ will be  semistable if $\deg F=\deg E/F$.

1st case: For each $p\in \Sigma$ the two eigenvalues of the local 
monodromy around $p$ coincide. Following Deligne, one can choose 
$(E,\nabla)$ such that at each $p\in \Sigma$ the two residue eigenvalues 
coincide. Thus in particular, $\deg F = \deg E/F$, and $E$ is semistable.

2nd case: For some $p\in \Sigma$ the two eigenvalues of the local monodromy
around $p$ differ. Let $(E,\nabla)$ be again a  Deligne
extension. Then the space 
$E\otimes \C(p)$ splits into two one-dimensional eigenspaces
$F\otimes \C(p)$ and $(E/F)\otimes \C(p)$ of the residue endomorphism
$\res_p(\nabla)$.
One can apply Gabber's construction \cite{EV} (section 1)
to either one of these eigenspaces and increase by one either the degree
of $F$ or that of $E/F$. Repeating this one can obtain bundles
$E'\supset F'$ with logarithmic connections such that ${\rm deg}
F' = {\rm deg} (E'/F')$.
Then $E'$ is semistable.
\end{proof}

\begin{rem}
If $\rho :\pi_1(X\setminus \Sigma)\to GL(2,\C)$ is reducible and for any $p\in \Sigma$
the local monodromy around $p$ has only one Jordan block then any semistable
bundle $(E,\nabla)$ with logarithmic connection which realizes $\rho$ 
is at each point $p\in \Sigma$ a Deligne extension by theorem \ref{ss}.
It satisfies $\deg E = 2\deg F \in 2\Z$. 

Examples of such representations are given in (2.1) ($n=2$) and
in (2.4). Or simply take $0\neq \alpha\in H^0(X, \omega)$ and
the connection 
\begin{gather*}
(\sO\oplus \sO, d+ \Big(\begin{array}{ll}
0 & \alpha\\
0 & 0\end{array}\Big))
\end{gather*}
if the genus is $\ge 1$.
\end{rem}

\section{Some three-dimensional representations}

Bolibruch's first class of representations
$$
\rho : \pi_1(X\setminus \Sigma )\to GL(n,\C)
$$
for $X=\P^1_\C,$ $\Sigma\subset X$ finite, which cannot be realized by
a semistable bundle with a logarithmic connection (with poles only in $\Sigma$)
has the following properties \cite{AB} (ch. 2):
\begin{list}{}{}
\item[(i)]
$\rho$ is three-dimensional and reducible with a one-dimensional 
subrepresentation $\rho'$.
\item[(ii)]
For each $p_i\in \Sigma$ the local monodromy of $\rho$ around $p_i$ has only
one eigenvalue $\lambda_i$ and one Jordan block.
\item[(iii)]
If $(E'', \nabla'')$ realizes $\rho'' := \rho/\rho'$ and if it is a Deligne
extension at each point $p\in \Sigma$ then $E''$ is not semistable.
\end{list}

By theorem \ref{ss} it is obvious that $\rho$ with (i) -- (iii) cannot
be realized by a semistable bundle with logarithmic connection
(with poles only in $\Sigma$).

\begin{rem}
(iii) follows if one knows a single bundle $(E'',\nabla'')$ with logarithmic
connection 
which realizes $\rho/\rho'$, which is a Deligne extension at each 
$p\in \Sigma$,
and which is not semistable. Then any other bundle which realizes
$\rho/\rho'$ and which is a Deligne extension at each $p\in \Sigma$
is obtained from $E''$ by tensoring with a suitable line bundle.
\end{rem}

\begin{rem}
If $\rho''$ is a two-dimensional representation with (ii) and
$\prod_i \lambda_i=1$ and $|\Sigma|\geq 2$ then one can construct easily a 
three-dimensional representation $\rho$ with (i) and (ii) and
$\rho/\rho' = \rho ''$.

We use the notations of section 2. Let $\rho''$ be given by 
$2\times 2$-matrices $A_1'', B_1'',...,A_g'',B_g''$ and 
$C_1'',...,C_m''$. Define
\begin{eqnarray*}
&& A_i := \left( \begin{array}{cc} 
1 & \begin{array}{cc} 0&0\end{array} \\
\begin{array}{c}0\\0\end{array} & A_i'' 
\end{array}\right) \ ,\ 
B_i := \left( \begin{array}{cc} 
1 & \begin{array}{cc} 0&0\end{array} \\
\begin{array}{c}0\\0\end{array} & B_i'' 
\end{array}\right) \ ,\\
&& C_i := \left( \begin{array}{cc} 
\lambda_i & \begin{array}{cc} \gamma_{i1}&\gamma_{i2}\end{array} \\
\begin{array}{c}0\\0\end{array} & C_i'' 
\end{array}\right)
\end{eqnarray*}
for suitable $\gamma_{ij}$ such that
$$
A_1B_1A_1^{-1}B_1^{-1}\cdot ...\cdot A_gB_gA_g^{-1}B_g^{-1}\cdot
C_1...C_m = \id
$$
holds. $\gamma_{1j},...,\gamma_{m-1j}$ can be chosen freely.
$\gamma_{m1}$ and $\gamma_{m2}$ are given by two linear functions
in $\gamma_{1j},...,\gamma_{m-1j}$ such that for each $i=1,...,m-1$
the linear parts in $\gamma_{i1},\gamma_{i2}$ of the two functions
are together invertible. For generic solutions $\gamma_{1j},...,\gamma_{mj}$
the matrices $C_i$ have only one Jordan block.
\end{rem}

Bolibruch proved (iii) for his examples by quite involved explicit 
calculations.
Other examples, for higher genus curves $X$ can be obtained 
as follows.

Let 
$f:Z\to X$ be a proper semistable, nonisotrivial  family of elliptic curves 
over a curve $X$.
Let $Y \subset Z$ be the union of the bad fibers. The Gau{\ss}-Manin bundle
$$
R^1f_{*}\Omega^{\bullet}_{Z/X}(\log Y)
$$
on $X$ has rank 2, and the Gauss-Manin connection $\nabla$ on it has 
logarithmic poles with milpotent residues 
along $ \Sigma \subset f(Y)$ ($\Sigma$ might be smaller, due to
bad fibers of $f$ inducing good fibers for the Jacobian family).
It contains the positive subbundle $f_*\omega_{Z/X}$,
thus is instable.
Once such an $f$ is chosen, one obtains other ones by
considering the pullback family over any covering of $X$, \'etale on
$\Sigma$. In particular, one can make the genus of $X$ arbitrarily high.

\end{document}